\documentclass[a4paper,12pt]{amsart}
\usepackage{amsmath}
\usepackage{amssymb}
\usepackage{amsthm}
\usepackage{graphicx,color,wrapfig}

\theoremstyle{plain} 
\newtheorem{lemma}{Lemma}
\newtheorem{Definition}[lemma]{Definition}
\newtheorem{prop}[lemma]{Proposition} 
\newtheorem{cor}[lemma]{Corollary}
\newtheorem{rem}[lemma]{Remark}

\newenvironment{Proof}[1][Proof]{
    \par
    \topsep 6pt plus 6pt
    \trivlist
    \item[\hskip\labelsep\bfseries #1.]\ignorespaces}%
    {\qed\endtrivlist
}

\newcommand{\R}{\mathbb{R}}

\newcommand{\M}{\mathbb{M}}

\begin{document}

\author[John Andersson, Norayr Matevosyan, Hayk 
Mikayelyan]{John Andersson, Norayr Matevosyan, Hayk Mikayelyan$^\dagger$}
\title[The tangential touch between free and  
fixed boundaries]{On the tangential touch between the free and the 
fixed boundaries for the two-phase obstacle-like problem}



\thanks{$\dagger$  Supported by
Deutsche Forschungsgemeinschaft and Freistaat Sachsen}

\thanks{Second and third authors thank G\"oran Gustafsson Foundation for 
visiting appointments to KTH, Stockholm}

\address{John Andersson\\
Institutionen f\"or Matematik \\
Kungliga Tekniska H\"ogskolan\\
100 44 Stockholm \\
Sweden}
\email{johnan@e.kth.se}

\address{Norayr Matevosyan\\
Johann Radon Institut f\"ur Angewandte Mathematik \\ 
\"Osterreichische Akademie der Wissenschaften \\
Altenbergerstrasse 69 \\
4040 Linz \\ 
Austria  }
\email{norayr.matevosyan@oeaw.ac.at}

\address{Hayk Mikayelyan \\
Mathematisches Institut \\
Universit\"at Leipzig\\\,\,\,
Augustusplatz 10/11\\
04109 Leipzig \\
Germany}
\email{hayk@math.uni-leipzig.de}

\begin{abstract}
\noindent In this paper we consider the following two-phase 
obstacle-problem-like equation in the unit half-ball
 \begin{equation*}
\Delta u = \lambda_{+}\chi_{\{u>0\}}-\lambda_{-}\chi_{\{u<0\}} ,\,\,
\lambda_\pm>0.
\end{equation*}
We prove that the free boundary touches the
fixed one in (uniformly) tangential fashion if the boundary data $f$
and its first and second derivatives vanish at the touch-point.
\end{abstract}

\maketitle

\section{Introduction}

\subsection{The Problem}

The following two-phase analogue of classical obstacle problem was suggested
by G. S. Weiss in \cite{W2} and then considered by N.N. Uraltseva in \cite{U} and
H. Shahgholian, N.N. Uraltseva and G.S. Weiss in \cite{SUW}. 
Study properties of a weak solution $u\in W^{1,2}(D)$ of
\begin{equation}
\Delta u = \lambda_{+}\chi_{\{u>0\}}-\lambda_{-}\chi_{\{u<0\}},
\label{main} 
\end{equation}
in the domain $D$, such that $u-f\in W_0^{1,2}(D)$ for a given $f\in W^{1,2}(D)$.
In our paper we always assume $\lambda_\pm>0$, and we consider
the cases $D=B_1$ and $D=B_1^+$, as well as the case of the
so-called global solutions $D=\R^n_+$. 

Obviously (\ref{main}) is the Euler-Lagrange equation of the energy functional 
$$
J(u)=\int_D |\nabla u|^2+2\lambda_+\max(u,0) +2\lambda_-\max(-u,0) dx.
$$
Note that if the boundary data $f$ is non-negative (non-positive) then 
the solution $u$ is so, too, and we arrive at classical obstacle problem
(see \cite{C}). In the two-phase case we do not have the property that
the gradient vanishes on the free boundary $\Gamma_u$ (see Section \ref{Notat}
for definition), as it was in the classical case; this causes difficulties.

We consider the following problem: Let $u$ be a weak solution of 
(\ref{main}) in $B_1^+$, $0\in\overline{\Gamma}_u$, 
$f:=u|_\Pi\in C^{2,Dini}(B_1\cap\Pi)$ and 
\begin{equation}
\label{cond}
f(0)=|\nabla f(0)|=|D^2 f(0)|=0.
\end{equation}
Then we prove that 
the free boundary of $u$ approaches the fixed one at $0$ tangentially.
Under some growth assumptions we prove that this approach is 
uniform (Theorem B) and we show the necessity of this assumption 
with an example.
From (\ref{cond}) it obviously follows that 
$\frac{|f(x')|}{|x'|^2}\leq \omega(|x'|)$ 
for some Dini modulus of continuity $\omega$,
i.e., the blow-up of $f$ is zero
$$
f_r(x'):=\frac{f(rx')}{r^2}\to 0 \,\,\,\text{as $r\to 0$}.
$$

Let us recall the definition of $C^{2,Dini}(B_1\cap\Pi)$; these are
functions from $C^2(B_1\cap\Pi)$ such that 
$$
|D^2 f(x)-D^2 f(y)| \leq \omega(|x-y|),
$$
where $D^2 f$ is the Hessian of $f$ and $\omega$ 
is a Dini modulus of continuity, i.e.,
$$
\int_0^1\frac{\omega(s)}{s}ds < \infty.
$$

\vspace{5mm}

\subsection{Notations}\label{Notat} 
In the sequel we use following notations: 
\vspace{3mm}

\begin{tabbing}
$ \mathbf{R}_{+}^{n}$  \hspace{25mm}      \=  $\{x\in \mathbf{R}^{n}:x_{1}>0\}$\\
$ \mathbf{R}_{-}^{n} $  \>  $ \{x\in \mathbf{R}^{n}:x_{1}<0\} $, \\
$  B(z,r) $  \>  $ \{x\in \mathbf{R}^{n}:|x-z|<r\} $, \\
$ B_{r} $  \>  $ B(0,r) $, \\
$ B^{+}_r $  \>  $ \mathbf{R}_{+}^{n}\cap B_r $, \\
$ \Pi  $  \>  $ \{x\in \mathbf{R}^{n}:x_{1}=0\} $, \\
$ x' $  \>  $ (x_2,\dots,x_n) $, \\
$ K_\epsilon $  \>  $ \{x\in \mathbf{R}^{n}_+:x_{1}>\epsilon |x'|\} $, \\
$ \Vert \cdot \Vert _{\infty } $  \>  canonical norm, \\
$ e_{1},\ldots ,e_{n} $  \>  standard basis in $  \mathbf{R}^{n} $, \\
$ \nu ,\ e  $  \>  arbitrary unit vectors, \\
$ D_{\nu },\ D_{\nu e} $  \>  first and second directional derivatives, \\
$ v^{+},\ v^{-} $  \>  $ \max (v,0),\ \max (-v,0) $, \\
$ \chi _{D} $  \>  characteristic function of the set $ D $, \\
$ \partial D $  \> boundary of the set $ D $, \\
$ \Omega ^{+}_u $  \>  $ \left\{ x\in D :u\left( x\right)>0\right\} $, \\
$ \Omega ^{-}_u $  \>  $ \left\{ x\in D :u\left( x\right)<0\right\} $, \\
$ \Lambda_u $  \>  $\left\{ x\in B_{1}^{+}:u\left( x\right)
=\left\vert \nabla u\left( x\right) \right\vert =0\right\}  $, \\
$ \Gamma_u $  \>  $ (\partial \Omega ^{+}_u\cup\partial 
\Omega ^{-}_u)\cap D, $ the free boundary, \\
$ \mathcal{P}(\dots) $  \>  see Definition \ref{defdef}. 
\end{tabbing}

\newpage

\subsection{``Typical'' examples}

We show here with some examples how the situation near a touch
point between free and the fixed boundaries can look like. 

\begin{figure}
\begin{center}
\input{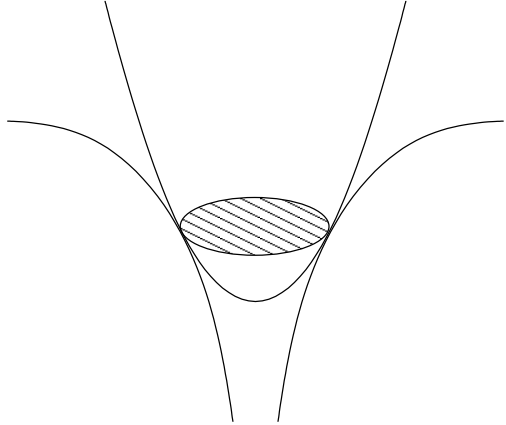}
\caption{``Ball solutions''}
\label{ballsol}
\end{center}
\end{figure}

Fix the ball $B_{R}$ and 
consider the function $\dfrac{\lambda}{2n}|x|^2$, $\lambda>0$. Let us take
the radial fundamental solution of Laplace equation $U$ multiplied with a 
constant $C_R$, such that $C_R \partial_r U(R)=\frac{\lambda}{n} R$. 
Then for some constant $C$ the
function 
$$
V(x)=\frac{\lambda}{2n}|x|^2-C_R U(|x|)+C
$$  
is non-negative in $\R^n$, 
$\Delta V=\lambda - C_R\delta_0 $ and $V=|\nabla V|=0$ on $\partial B_R$
(Figure \ref{ballsol}).

\begin{figure}
\begin{center}
\input{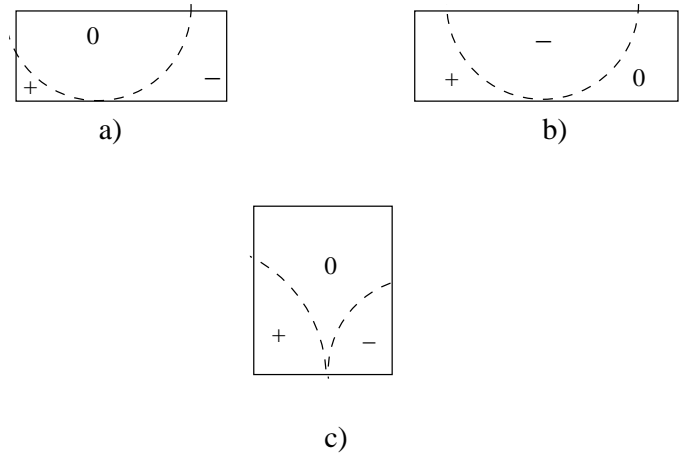}
\caption{``Typical'' examples}
\label{typex}
\end{center}
\end{figure}

Now we see how we can construct solutions of 
(\ref{main}) in $B_R\backslash \{0\}$
and in $\R^n\backslash B_R$. So for instance in $\R^2$ we can illustrate
some solutions considered in rectangles (Figure \ref{typex}). 
The dashed curves denote free boundaries, $\pm$ denote regions 
$\Omega_u^\pm$ and $0$ the region $\Lambda_u$. 
Figure \ref{typex} a) shows that the case
when the solution does not have quadratic growth 
near the touch point is possible. In this case the blow-up of the 
solution is zero.
Figure \ref{typex} b) shows that even if we have non-negative boundary data near
touch point, the blow-up still can be negative and Figure \ref{typex} c) shows that
the condition (\ref{cond}) is essential for the tangential touch. 

Let us take boundary data $f$ on $\partial B_1^+$ to be odd-symmetric 
with respect to $x_2$. Then the solution $u$ will be odd-symmetric, too. 
From our results (see Section \ref{MAINRES}) follows 
that the set $\Lambda_u$ is large near 0, where the free boundary
touches the fixed boundary, as it is
the case in the example from the Figure \ref{typex} a). So we do not have 
orthogonal touch, as it might be expected. The similar argument works 
also in higher dimensions for every plane-symmetric domain.

\newpage
\section{Main results}\label{MAINRES}

In this section we state two theorems. The first one says that
if the boundary data of a solution of (\ref{main}) satisfies condition 
(\ref{cond}), then the free boundary can approach the fixed one 
only tangentially. In the second theorem we assert that this approach is
uniform for a certain class of functions.

\vspace{5mm}

\noindent \textbf{Theorem A}. \textsl{\ \ Let $u$ be a solution of (\ref{main})
in $B_1^+$ with boundary data $f$ on $\Pi$, condition (\ref{cond})
is satisfied and $0\in \overline{\Gamma}_u $.
Then the free boundary approaches $\Pi$ at the point $0$ tangentially.
}

\vspace{3mm}

\begin{cor} 
\label{corB}
Let $u$ be as in Theorem A, then 
one of the following expressions is true
\begin{equation*}
\dfrac{|\Omega^+_{u}\cap B_r^+|}{|B_r^+|}\to 1,\,\,
\dfrac{|\Omega^-_{u}\cap B_r^+|}{|B_r^+|}\to 1,\,\,
\dfrac{|\Lambda_{u}\cap B_r^+|}{|B_r^+|}\to 1\,\, \text{as}\,\,r\to 0.
\end{equation*}
Moreover, the first two cases are possible only if the condition
(\ref{paym2}), see below, is satisfied for some $c_0$ and $r_0$, 
and the third case if it fails to hold (see Lemma \ref{Lemnondeg}).
\end{cor}

\vspace{3mm}

\begin{Definition}
\label{defdef}
Let $\omega$ be a Dini modulus of continuity 
and $M$, $c_0$, $r_0$ be positive constants.
We define 
$\mathcal{P}(M,R,c_0,r_0)$ to be the 
class of solutions $u$ of (\ref{main}) in $B_1^+$,
$\|u\|_{L^\infty(B_1^+)}\leq M$, $0\in \overline{\Gamma}_u$ such that
the boundary data 
$f=u|_\Pi\in C^{2,Dini}(B_1\cap\Pi)$ satisfies condition (\ref{cond}) and 
\begin{equation}
\label{condstar}
\|f\|_{C^2(\overline{B}_1\cap\Pi)}\leq R,\,\,\,
\int_0^1\frac{\omega(s)}{s}ds\leq R.
\end{equation}
We assume, further
\begin{equation}
\sup_{B_{r}^+}|u|\geq c_0 r^2,\,\,\text{for} \,\,0<r<r_0,
\label{paym2} 
\end{equation}
for all $u\in\mathcal{P}(M,R,c_0,r_0)$.
\end{Definition}

\begin{figure}
\begin{center}
\input{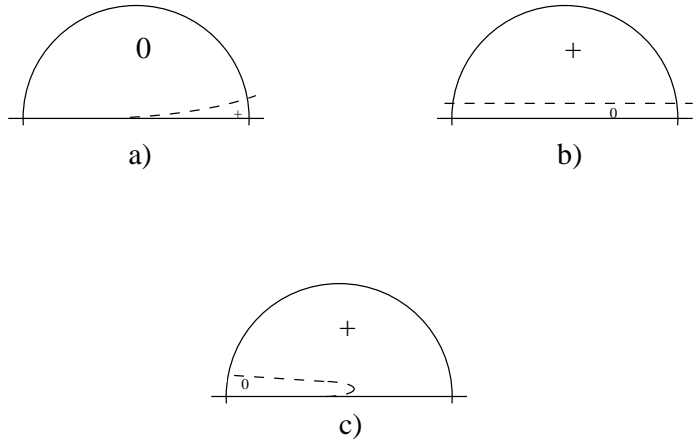}
\caption{Non-uniform approach}
\label{nonunif}
\end{center}
\end{figure}
\vspace{1mm}

\begin{rem}
If $u$ solve (\ref{main}) in $B_1^+$, 
$0\in\overline{\Gamma}_u$ and $u|_\Pi\equiv 0$, then
condition (\ref{paym2}) is fulfilled with the constant
$c_0=max(\lambda_\pm)C$, for $0<r<1$, where $C$ is a 
dimension dependent constant (see Lemma \ref{qaraj}
and Corollary \ref{corqor}).
\end{rem}

\vspace{5mm}

\noindent \textbf{Theorem B}. \textsl{\ \ 
There exists a modulus of continuity 
$\sigma(r)$ and ${\tilde r}>0$ such that if
$u\in \mathcal{P}(M,R,c_0,r_0)$ then
$$
\Gamma_u\cap B_{{\tilde r}}\subset \{x:x_1<|x'|\sigma(|x'|)\}.
$$
In other words the free boundary of the functions
from $\mathcal{P}$ approaches $\Pi$ at the 
point $0$ uniformly tangentially.
}

\vspace{2mm}
Here $\sigma$ and ${\tilde r}$ depend on the dimension, $\lambda_\pm$,
$c_0$, $r_0$, $M$ and  $R$.

\begin{rem}
Since the main tool we use proving Theorems A 
and B is the blow-up argument, these results 
could be generalized for domains with smooth 
enough boundary. 
\end{rem}

The following example shows that we do not 
have uniform tangential approach
when the set $\Lambda_u$ is large near the 
touch point, i.e., the condition
(\ref{paym2}) fails to hold (see Lemma \ref{Lemnondeg}). 
In the example
the boundary data is positive, so we treat here
the classical obstacle problem in $\R^2$.
First let us take small positive boundary data on the right
side from touch point (zero) to get a tiny 
positivity set as it is shown on the  
Figure \ref{nonunif} a). This can be done using
the ``ball solutions'' above. Next consider the function 
$u(x)=\frac{\lambda_+}{2}((x_1-\epsilon)_+)^2$.
We will get it as a solution if we take its boundary 
data on $\partial B^+_1$ (Figure \ref{nonunif} b)).
Consider now the boundary data which is the sum of the boundary data
of the previous two examples. 
It will look like it is shown on the Figure \ref{nonunif} c).
So we see that when the `tiny' positivity set in Figure \ref{nonunif} a) 
becomes smaller and $\epsilon$ in Figure \ref{nonunif} b)
tends to zero, then for any ${\tilde r}$ and $\sigma$ 
in the Figure \ref{nonunif} c) we get free boundary points 
in $ B_{\tilde r}\cap\{x:x_1>|x'|\sigma(|x'|)\}$ , while the boundary 
data on $\Pi$ remains bounded.

\vspace{3mm}

\begin{rem}
In order to get uniform tangential touch for a class of
solutions we impose condition (\ref{paym2}). This condition,
however, can be replaced by the following one, which is 
considered in \cite{KKS} for a different problem;
\begin{equation*}
\dfrac{|\Omega_{u}^+\cap B_r^+|}{|B_r^+|}\geq c_0>0,\,\, 
\text{for}\,\,r<r_0.
\end{equation*}
From Lemma \ref{Lemnondeg} and Corollary \ref{corB} it
follows that both conditions are equivalent in our case.
\end{rem}

\vspace{5mm}

\section{Technicalities}
\subsection{Non-degeneracy}
In this section we 
introduce some (modified) results from \cite{W2}, \cite{U}
and \cite{SUW} as well as prove growth estimates at the 
boundary (Lemmas \ref{Lemnondeg} and \ref{AndLem}).

\begin{lemma}
\label{qaraj}
Let $u$ solve (\ref{main}) in $B_1$. There exists a dimension dependent
constant $C$ such that
$\|f^\pm\|_\infty<\lambda_\pm C$ implies 
$0\notin \overline{\Omega}_u^\pm$.
\end{lemma}
\begin{Proof}
Consider the ``+''-case.
Due to the comparison principle a similar argument is true (and well-known)
for obstacle problems, i.e. it is known in our case if the boundary data 
$f$ is non-negative or non-positive. 
Let us now consider the related (one-phase) obstacle problem
in $B_1$ with boundary data $f^+$, denote its solution by $v$. 
It is enough to show 
that $\Omega^+_u\subset\Omega^+_v$. Consider the function
$w=u-v$ in $\Omega^+_u$. We have $w|_{\partial\Omega^+_u}\leq 0$
and $\Delta w\geq 0$ in $ \Omega^+_u$, hence $u-v\leq 0$ and we are done.
\end{Proof}

\begin{cor}
\label{corqor}
Let $u$ be the solution of (\ref{main}), $x_0\in\overline{\Omega}^\pm$
and $B_{r_0}(x_0)\subset D$. Then 
\begin{equation}
\label{nondeg-1} 
\sup_{\partial B_{r}(x_0)} u^\pm\geq \lambda_\pm C r^2,\,\,\text{for} \,\,r<r_0.
\end{equation}
Here the constant $C$ is the same as in the previous lemma.

\noindent In other words if
$x_0\notin\text{int}\,\, \Lambda_{u}$ and
$B_{r_0}(x_0)\subset D$, then 
\begin{equation}\label{nondeg} 
\sup_{\partial B_{r}(x_0)}|u|\geq\min(\lambda_\pm) C r^2,\,\,\text{for} \,\,r<r_0.
\end{equation}
\end{cor}
\begin{Proof}
Let us restrict the function $u$ to $B_{r}(x_0)$ and scale 
\begin{equation*}
u_r(x)=\dfrac{u(rx+x_0)}{r^2}.
\end{equation*}
$u_r$ is then a solution of (\ref{main}) in $B_1$ with boundary data
$u_r|_{\partial B_1}$. Since $0\in\overline{\Omega}^+$ we must have
\begin{equation*}
\sup_{\partial B_1}u_r\geq \lambda_+ C, 
\end{equation*}
which in turn implies (\ref{nondeg-1}).
\end{Proof}

\begin{lemma}
\label{Lemnondeg} 
Let $u$ be the solution of (\ref{main}) in $B_1^+$ and suppose
for given constants $c_0$, $r_0$, we have
\begin{equation*}
\dfrac{|\Omega_{u}^+\cap B_r^+|}{|B_r^+|}\geq c_0>0,\,\, \text{for}\,\,r<r_0.
\end{equation*}
Then there exists a constant $c$ depending on, $c_0$, $\lambda_\pm$ and
the dimension such that 
\begin{equation}
\label{NONdeg} 
\sup_{B_{r}^+}|u|\geq c r^2,\,\,\text{for} \,\,r<r_0.
\end{equation}
\end{lemma}
The same is also true for $\Omega^-$.
\begin{Proof}
Let us denote ${\tilde B}^+_r=B_r^+\cap\{x_1>\epsilon r\}$. We can fix
an $\epsilon>0$ such that 
\begin{equation*}
\dfrac{|\Omega_{u}^+\cap {\tilde B}^+_r|}{|{\tilde B}^+_r|}\geq \frac{c_0}{2},
\,\, \text{for}\,\,r<r_0.
\end{equation*}
Hence for each $r>0$ there exists 
$x_r\in \Omega_{u}^+\cap{\tilde B}^+_{\frac{r}{2}}$.
Applying previous corollary to the ball $B_{d_r}(x_r)$, 
where $d_r$ is the $e_1$
component of $x_r$ we get that
\begin{equation*}
\sup_{B_{r}^+}|u|\geq \sup_{B_{d_r}(x_r)}|u|
\geq \lambda_+ C d_r^2\geq \epsilon\lambda_+\frac{ C }{4}r^2.
\end{equation*}
We proved the lemma with $c=\epsilon\lambda_+\frac{  C }{4} $.
\end{Proof}
In the proof of the next lemma we use the 
technique from \cite{A} (Lemma 5), see also \cite{CKS}.
A similar estimate in the interior was 
proved by Uraltseva in \cite{U}.

\begin{lemma}
\label{AndLem}
Let $u$ solve (\ref{main}) in $B_1^+$, $\|u\|\leq M$ 
and assume its boundary 
data $f=u|_\Pi$ and the Dini modulus of continuity $\omega$ 
satisfy conditions (\ref{cond}) and (\ref{condstar}). 
Then there exists a 
constant $C=C(M,R)$ such that 
$$
\sup_{B_r^+}|u-D_{e_1}u(0)x_1|\leq C r^2,\,\,\,0<r<\frac{1}{2}.
$$
\end{lemma}

\begin{Proof}
Let us denote by 
$$
S_j(u):=\sup_{B_{2^{-j}}^+}|u-D_{e_1}u(0)x_1|
$$
and $\M(u):=\{j:S_j(u)\leq 4 S_{j+1}(u)\}$.
We want to show that $S_j(u)\leq C 2^{-2j}$.
First let us show this for all $j\in\M(u)$. The proof is done 
by contradictory argument: assume there exist a sequence $\{u_j\}$ of solutions of 
(\ref{main}) in $B_1^+$ such that
$$
S_{k_j}(u_j)\geq j 2^{-2 k_j},
$$
for some $k_j\in\M(u_j)$. Denoting by $w_j(x):=u_j(x)-D_{e_1}u_j(0)x_1$
and by 
$$
{\tilde w}_j(x):=\dfrac{w_j(2^{-k_j}x)}{S_{k_j+1}(u_j)},
$$
we get 
$$
\|\Delta{\tilde w}_j\|_\infty\leq \text{max}(\lambda_\pm)
\frac{2^{-2k_j}}{S_{k_j+1}(u_j)}\leq\text{max}(\lambda_\pm)
\frac{2^{-2k_j}}{\frac{1}{4}S_{k_j}(u_j)}\leq\text{max}(\lambda_\pm)
\frac{4}{j}\to 0.
$$
We also have 
\begin{equation}
\label{suptild1} 
\sup_{B_{\frac{1}{2}}^+}|{\tilde w}_j|=1.
\end{equation}
The condition
$(D_e f(x'))^\pm\leq |x'|\omega(|x'|) $,
for any unit vector $e\in \Pi$, implies that
\begin{equation}
\label{GreenEst}
\sup_{B^+_r}|D_e w_j|\leq Cr,
\end{equation}
where $C$ depends on $M$ and $R$.
To see this we should consider harmonic functions
$v^\pm_j$ in $B^+_{\frac{1}{2}}$ with the same boundary 
data as $(D_e w_j)^\pm$. Inequality 
(\ref{GreenEst}) then follows from
the subharmonicity of $(D_e w_j)^\pm$ (see \cite{U}) and
standard estimates on Green's function for the half-ball
(see \cite{Wi}). From (\ref{GreenEst}) we have
\begin{equation}
\label{Crj} 
\sup_{B^+_r}|D_e {\tilde w}_j|\leq \frac{4Cr}{j},
\end{equation}
A subsequence of ${\tilde w}_j$ converges in 
$C^1(B^+_{\frac{1}{2}})$ to
a harmonic function $u_0$. Due to (\ref{Crj}) we get
$D_e u_0=0$ for all $e\in\Pi$, thus $u_0=ax_1$. On
the other hand $D_{e_1}{\tilde w}_j(0)=0$
and by $C^1$-convergence (up to $\Pi$) the same holds for 
$u_0$. Hence $u_0\equiv 0$ which contradicts (\ref{suptild1}).

Next let us show that $S_j(u)\leq 4C 2^{-2j}$ for all $j$.
Suppose $j$ is the first integer for which the inequality fails
to hold, then
$$
S_{j-1}(u)\leq 4 C 2^{-2(j-1)}\leq 4 S_j(u),
$$
i.e. $j-1\in\M(u)$ and 
$$
S_j(u)\leq S_{j-1}(u)\leq C 2^{-2(j-1)}=4C2^{-2j},
$$
a contradiction.
\end{Proof}

\subsection{Monotonicity formulae}

Here we introduce two monotonicity formulae in the
following two lemmas, which
play crucial role in our proofs. The first one was 
presented by H. W. Alt, L. A. Caffarelli and A. Friedman in
\cite{ACF} and was developed then in \cite{CKS}. 
The second one is due to G. Weiss \cite{W1}, \cite{SUW}. 
In \cite{A} Andersson adapted it to the half-space case and
our representation is analogous. See also \cite{M} for the
formula in parabolic case.
\begin{lemma}
{\bf (ACF monotonicity formula)}

Let $h_1$, $h_2$ be two non-negative continuous subsolutions of
$\Delta u=0$ in $B_R$. Assume further that $h_1 h_2=0$ and
that $h_1(0)=h_2(0)=0$. Then the following function is 
non-decreasing
in $r\in(0,R)$
\begin{equation}
\varphi(r)=\dfrac{1}{r^4}
\left( 
\int_{B_r}\dfrac{|\nabla h_1|^2 dx}{|x|^{n-2}}
\right)
\left( 
\int_{B_r}\dfrac{|\nabla h_2|^2 dx}{|x|^{n-2}}
\right).
\end{equation}
More exactly, if any of the sets $\text{spt}(h_j)\cap\partial B_r$
digresses from a spherical cap
by a positive area, then either $\varphi'(r)>0$ or $\varphi(r)=0$.
\end{lemma}

\begin{lemma}
{\bf (Weiss' monotonicity formula)}

Let $u$ solve (\ref{main}) in $B_R^+$ and $u|_{\Pi\cap B_R}=0$.
Then the function
\begin{equation}
\Phi(r)=r^{-n-2}\int_{B_r\cap\R^n_+} (|\nabla u|^2+2\lambda_+u^+ +2
\lambda_-u^-)-r^{-n-3}\int_{\partial B_r\cap \R^n_+} 2u^2d{\mathcal H}^{n-1}
\end{equation}
is non-decreasing for $r\in (0,R)$. Moreover, if 
$\Phi(\rho)=\Phi(\sigma)$ for any $0<\rho<\sigma<R$ then 
$\Phi$ is homogeneous of degree two in 
$(B_\sigma\backslash B_\rho)\cap \R^n_+$.
\end{lemma}
The proof is analogous to the proof of the Lemma 1 in \cite{A}.

\section{Global solutions}

In this section we will classify all solutions of (\ref{main})
in the $\R^n_+$ with zero boundary data and quadratic growth.
We will see that only possible cases are 
\begin{equation}
u(x)=\pm \dfrac{\lambda_\pm}{2}((x_1-a)_+)^2,  \,\,\, a\geq 0\,\,\,
\text{or}\,\,\, u(x)=\pm \dfrac{\lambda_\pm}{2}x_1^2\pm\alpha x_1, \,\,\,
\alpha\geq 0.
\label{pmux}
\end{equation}
The proofs of next two lemmas adapt the proofs of analogous 
results from \cite{SU} for our case.

Let us first prove that $u$ is two-dimensional.
\begin{lemma}
\label{twodim} 
Let $u$ solve (\ref{main}) in $\R^n_+$
with boundary data $u|_\Pi=0$.
Then the function $u$ is two-dimensional, 
i.e., in some system of coordinates
$$
u(x)=u(x_1,x_2),
$$
where the $e_1$ direction is normal to $\Pi$.
\end{lemma}

\begin{Proof}
Let us take any $e$ orthogonal to $e_1$ and consider functions $(D_e u)^\pm$.
In \cite{U} Uraltseva proved that these functions are subharmonic. Note that they
will remain so if we extend them by zero to $\R_-^n$. Now we can apply 
ACF monotonicity formula to $(D_e u)^\pm$. For $r<s$ we have
$$
\varphi(r,D_e u)\leq\varphi(s,D_e u)\leq 
\lim_{s\to\infty}\varphi(s,D_e u)=:C_e.
$$
In \cite{U} is shown that the second derivatives of $u$ are bounded, thus
we can find a sequence $u_{r_j}=\frac{u(r_j x)}{r_j^2}\to u_\infty$,
uniformly on compact subsets and in 
$(W^{2,p}_{loc}\cap C^{1,\alpha}_{loc})(\R^n_+\cup\Pi)$, for any 
$1<p<\infty$ and $0<\alpha<1$.
We have now
$$
C_e=\lim_{r_j\to\infty}\varphi(sr_j,D_e u)=
\lim_{r_j\to\infty}\varphi(s,D_e u_{r_j})=
\varphi(s,D_e u_\infty),\,\,\,\forall s>0.
$$
From $\{x_1<0\}\subset \{D_e u=0\}$ follows that
$\varphi(r,D_e u_\infty)\equiv 0$ or $\varphi ' (r,D_e u_\infty)>0$
for all $r>0$, thus $C_e=0$ and we get $D_e u \geq 0$ or 
$D_e u \leq 0$.  

For $e_2\in \Pi$ assume $D_{e_2}u\geq 0$ and let
$e_3\in \Pi$ be orthogonal to $e_2$.
Consider unit vector $e(\phi)=\cos\phi \,e_2+\sin\phi \,e_3\in\Pi$, 
$\phi\in[0,\pi]$. From the $C^1$-continuity we have that
the sets $\{\phi:\Omega^\pm_{D_{e(\phi)}u}\not=\emptyset\}$ 
are relatively open in $[0,\pi]$. On the other hand they are both non-empty
and have empty intersection; this means that 
there exists $\phi_0\in(0,\pi)$ such that
$D_{e(\phi_0)}u\equiv 0$. Rotating coordinates we can get
$D_{e_2}u\geq 0$ and $D_{e_3}u\equiv 0$. Repeating this
with $e_k$, $k=4,\ldots ,n$, we get that $u$ is two dimensional.
\end{Proof}

We prove now the main result of this section 
under the assumption of homogeneity.

\begin{prop}
\label{homogsol} 
Let $u$ be homogeneous of degree two and solve (\ref{main}) in $\R^n_+$
with boundary data $u|_\Pi=0$.
Then either $u(x)=\frac{\lambda_+}{2}x_1^2$ or 
$u(x)=-\frac{\lambda_-}{2}x_1^2$.
\end{prop}
\begin{Proof}
We can consider only two-dimensional functions $u$. So let us rewrite
$u$ in radial coordinates as
$$
u(x)=u(r,\theta)=r^2 \phi(\theta),\,\,\, r\in [0,\infty),
\,\,\theta\in[0,\pi].
$$
Then we get the ODE
$$
\phi''+4\phi=\lambda_+\chi_{\{\phi>0\}}-\lambda_-\chi_{\{\phi<0\}}
$$
in the interval $[0,\pi]$ with boundary data $\phi(0)=\phi(\pi)=0$.
It can be checked that the only solutions of this
ODE are $\phi(\theta)=\pm\frac{\lambda_\pm}{2}\text{sin}^2\theta$.
\end{Proof}

\begin{lemma}
Let $u$ solve (\ref{main}) in $\R^n_+$
with boundary data $u|_\Pi=0$ and be quadratically bounded at infinity.
Then $u$ is one of the representations in (\ref{pmux}).
\label{lemglobsol} 
\end{lemma}
\begin{Proof}
If the function $u$ is non-negative or non-positive, then the 
result we want to prove follows from Theorem B in \cite{SU}. So let us show
that $u$ does not change the sign. We do this by contradiction; assume
that $u^\pm$ are both non-trivial. 

Consider the shrink down of $u$; ${\tilde u}:=\lim_{j\to\infty}u_j$,
where $u_j(x)=\frac{u(r_j x)}{r_j^2}$, $r_j\to\infty$. 
It is homogeneous of degree two. To see this
we need to use Weiss' monotonicity formula
$$
\Phi(s,{\tilde u})=\lim_{j\to\infty}\Phi(s , u_j)=
\lim_{j\to\infty}\Phi(s r_j, u)=\Phi(\infty,u).
$$
Thus ${\tilde u}$ equals to one of
$\pm \frac{\lambda_{\pm}}{2}x_1^2$ by Proposition \ref{homogsol} above. 
Assume for definiteness that we have the ``+''-sign.

This means that for any $\delta >0$ there exists $R_{\delta}$ such that 
\begin{equation}
\label{eq:sptcon}
\Omega_u^-\setminus B^+_{R_{\delta}} \subset \{ x;\; x_1<\delta |x_2| \}.
\end{equation}

Let us now take the barrier function
$$
U(x_1,x_2)= x_1^4 + x_2^4 - 6 x_1^2 x_2^2 + C.
$$ 
For large enough $C$ we have $\Omega_u^-\Subset \Omega^+_U$.
Since $u$ is quadratically bounded we get from the comparison 
principle that
$u^-(x)\leq\varepsilon U(x)$ for any $\varepsilon>0$, thus
$\Omega_u^-=\emptyset$.
\end{Proof}

\vspace{5mm}

\section{Proofs}
\begin{Proof}[Proof of the Theorem A]
We consider here only the case when (\ref{paym2}) fails to hold.
From Lemma \ref{AndLem} follows that $D_{x_1}u(0)=0$ and
\begin{equation}
\label{quadrbddness} 
\sup_{B_{r}^+}|u|\leq c_0 r^2,\,\,\text{for} \,\,r<r_0.
\end{equation}
Now assume that we do not have tangential touch at $0$, i.e.,
there is an $\epsilon>0$ and a sequence $x^j\in K_\epsilon\cap\Gamma_u$,
$x^j\to 0$. Repeating the proof of Lemma \ref{Lemnondeg}
we obtain
\begin{equation}
\sup_{B_{2d_j}^+}|u|\geq C d_j^2,\,\,\text{for} \,\,r<r_1,
\label{zrotoliche} 
\end{equation}
where $d_j=|x^j|$.
Consider now the blow up sequence
$$
{\tilde u}_j(x)=\frac{u(2d_jx)}{4d^2_j},
$$
which is bounded by (\ref{quadrbddness}). 
Therefore there is a subsequence converging in $C^{1,\alpha}$
to a global solution $u_0$ with zero boundary data,
which is non-trivial (due to (\ref{zrotoliche})). 
As in the proof of Lemma \ref{lemglobsol}, using Weiss' monotonicity formula 
we get that $u_0$ is homogeneous of degree $2$, this implies that 
$u_0(x)=\pm\frac{\lambda_\pm}{2}x_1^2$
which contradicts the fact that $x^j\in K_\epsilon$.

\end{Proof}

\begin{Proof}[Proof of the Theorem B]
The proof is done by contradictory argument. Assume there exist an $\epsilon>0$, 
functions $u_j$ satisfying the conditions of the theorem and a sequence
$x^j\to 0$ such that $x^j\in K_\epsilon\cap\Gamma_{u_j}$.
Let us consider the blow-up sequence
$$
{\tilde u}_j(x)=\frac{u_j(d_jx)}{\sup_{B^+_{d_j}}|u_j|},
$$
where $d_j:=|x^j|$.
We have that
$$
\Delta{\tilde u}_j= \frac{d_j^2}{\sup_{B^+_{d_j}}|u_j|}\Delta u_j.
$$
Two cases are possible: either
\begin{equation}
\frac{d_j^2}{\sup_{B^+_{d_j}}|u_j|}
\to 0
\label{dzgtzroi} 
\end{equation}
for some subsequence or 
\begin{equation}
\label{CHdzgtzroi} 
\frac{d_j^2}{\sup_{B^+_{d_j}}|u_j|}
\nrightarrow 0
\end{equation}
for all subsequences.

Let us consider the first case.
From Lemma \ref{AndLem} it follows that
\begin{equation}
-C r^2 + |D_{e_1}u_j(0)|r \leq{{\sup_{B^+_r}}}|u_j|
\leq C r^2 + |D_{e_1}u_j(0)|r.
\label{gnahat} 
\end{equation}
This together with (\ref{dzgtzroi}) gives 
that $|D_{e_1}u_j(0)|d_j^{-1}\to\infty$,
so we can assume 
\begin{equation}
|D_{e_1}u_j(0)|>j d_j.
\label{mecejic} 
\end{equation}
From here, and (\ref{gnahat}) we obtain
$$
\left|\dfrac{{\sup_{B^+_{d_j}}|u_j|}}{d_j |D_{e_1}u_j(0)|}-1\right|
\leq \frac{C}{j}\to 0.
$$
We arrive at
\begin{equation}
{{\sup_{B^+_{r}}}}|{\tilde u}_j|=
\frac{{\sup_{B^+_{rd_j}}|u_j|}}{{\sup_{B^+_{d_j}}|u_j|}}\leq
\frac{Cr^2d_j^2+|D_{e_1} u_j(0)|rd_j}{{\sup_{B^+_{d_j}}|u_j|}}
\to r.
\label{lastone} 
\end{equation}
There is a subsequence of ${\tilde u}_j$ converging to a function $u_0$
in $C^{1,\alpha}$, that 
is harmonic in $\R_+^n$ (due to (\ref{dzgtzroi})), linearly bounded 
(due to (\ref{lastone}))
and has zero boundary data at $\Pi$. Extending $u_0$ by odd reflection
to $\R_-^n$ and using Liouville's theorem we get that $u_0(x)=D_{e_1} u_0(0) x_1$
which contradicts the fact of existence of zeros in $K_\epsilon$.

In the case (\ref{CHdzgtzroi})
we can without loss of generality assume 
\begin{equation}
\frac{d_j^2}{\sup_{B^+_{d_j}}|u_j|}
\to d>0.
\label{notdztzroi} 
\end{equation}
We have then that a subsequence of ${\tilde u}_j$ converges to a
function $u_0$ in $C^1$ and (\ref{notdztzroi}) implies that
$u_0$ is a global solution with $d\lambda_\pm$ 
instead of $\lambda_\pm$ and zero boundary data.
Condition (\ref{paym2}) and Lemma \ref{lemglobsol}
give us that $u_0$ is strictly positive or negative in $\R^n_+$, which
contradicts the fact that $x^j\in K_\epsilon\cap\Gamma_{u_j}$.
More precisely, functions ${\tilde u}_j$ vanish at 
${\tilde x}_j:=d_j^{-1}x_j\in K_\epsilon\cap\Gamma_{u_j}\cap \partial B_1$,
thus we can always choose the subsequence of ${\tilde u}_j$ in such a way that
the corresponding subsequence 
${\tilde x}_j\to x_0\in K_\epsilon\cap\Gamma_{u_j}\cap \partial B_1$
and then $u_0(x_0)=0$.
\end{Proof}

\subsection*{Acknowledgment}
The authors are grateful to Prof. H. Shahgholian for valuable 
discussions and his kind hospitality.

\end{document}